\documentclass[11pt,a4paper]{article}

\usepackage{url,amsmath,amssymb,latexsym,pstricks,mathrsfs,comment,amsthm,graphicx,tikz,tikz-cd,enumerate,accents,pgffor,cite,wrapfig,multicol,float,cases,bibspacing,geometry}

\allowdisplaybreaks

\newcommand{\nc}{\newcommand}
\nc{\rnc}{\renewcommand}

\rnc{\emptyset}{\varnothing}
\nc{\ms}{\medskip}
\nc{\Part}{\mathscr P}
\nc{\dia}{\mathrel\diamond}
\nc{\D}{\mathcal D}
\nc{\T}{\mathcal T}
\rnc{\P}{\mathcal P}
\nc{\direl}[6]{\textstyle{\big[{#1\atop #4}{#2\atop#5}{#3\atop#6}\big]}}
\nc{\cL}{\mathcal L}
\nc{\cR}{\mathcal R}
\rnc{\H}{\mathscr H}
\rnc{\L}{\mathscr L}
\nc{\R}{\mathscr R}
\nc{\J}{\mathscr J}
\nc{\bn}{\mathbf{n}} 
\nc{\bA}{\mathbf{A}} 
\nc{\bB}{\mathbf{B}} 
\nc{\B}{\mathcal B} 
\nc{\I}{\mathcal I} 
\rnc{\S}{\mathcal S} 
\nc{\id}{\operatorname{id}}
\nc{\al}{\alpha}
\nc{\be}{\beta}
\nc{\ga}{\gamma}
\nc{\de}{\delta}
\nc{\set}[2]{\{ {#1} : {#2} \}} 
\nc{\bigset}[2]{\big\{ {#1} : {#2} \big\}} 
\nc{\sub}{\subseteq}
\nc{\rank}{\operatorname{rank}}
\nc{\coker}{\operatorname{coker}}
\nc{\codef}{\operatorname{codef}}
\nc{\deff}{\operatorname{def}}
\nc{\dom}{\operatorname{dom}}
\nc{\codom}{\operatorname{codom}}
\nc{\sm}{\setminus}
\nc{\COMMA}{,\quad}
\nc{\AND}{\qquad\text{and}\qquad}
\nc{\itemit}[1]{\item[\emph{(#1)}]}
\nc{\pf}{\noindent{\bf Proof.}  }
\nc{\pfof}[1]{\noindent{\bf Proof of #1.}  }
\nc{\epf}{\hfill$\Box$\bigskip}
\nc{\epfres}{\hfill$\Box$}
\nc{\bit}{\vspace{-3 truemm}\begin{itemize}}
\nc{\eit}{\end{itemize}\vspace{-3 truemm}}
\nc{\lam}{\lambda}
\nc{\si}{\sigma}
\nc{\Si}{\Sigma}
\nc{\Om}{\Omega}
\nc{\la}{\langle}
\nc{\ra}{\rangle}
\nc{\epfreseq}{\tag*{$\Box$}}
\rnc{\iff}{\ \Leftrightarrow\ }

\begin{document}

\numberwithin{equation}{section}

\newtheorem{thm}[equation]{Theorem}
\newtheorem{lemma}[equation]{Lemma}
\newtheorem{cor}[equation]{Corollary}
\newtheorem{prop}[equation]{Proposition}
\newtheorem{conj}[equation]{Conjecture}

\theoremstyle{definition}

\newtheorem{rem}[equation]{Remark}
\newtheorem{defn}[equation]{Definition}
\newtheorem{eg}[equation]{Example}
\newtheorem{ass}[equation]{Assumption}
\newtheorem{prob}{Problem}

\title{Ranks of ideals in inverse semigroups of difunctional binary relations}
\author{
James East\footnote{Centre for Research in Mathematics, School of Computing, Engineering and Mathematics, Western Sydney University, Locked Bag 1797, Penrith NSW 2751, Australia. {\it Email:} {\tt j.east\,@\,westernsydney.edu.au}} \ and 
Alexei Vernitski\footnote{Department of Mathematical Sciences, University of Essex, Colchester,
United Kingdom. {\it Email:} {\tt asvern\,@\,essex.ac.uk}}
}
\date{}

\maketitle

\vspace{-0.5cm}

\begin{abstract}
The set $\D_n$ of all difunctional relations on an $n$ element set is an inverse semigroup under a variation of the usual composition operation.  
We solve an open problem of Kudryavtseva and Maltcev (2011), which asks:  What is the
%
rank (smallest size of a generating set) of $\D_n$?  Specifically, we show that the rank of~$\D_n$ is  $B(n)+n$, where $B(n)$ is the $n$th Bell number.  We also give the rank of an arbitrary ideal of~$\D_n$.  Although $\D_n$ bears many similarities with families such as the full transformation semigroups and symmetric inverse semigroups (all contain the symmetric group and have a chain of $\mathscr J$-classes), we note that the fast growth of $\rank(\D_n)$ as a function of $n$ is a property not shared with these other families.

{\it Keywords}: Semigroups, binary relations, ideals, generators, rank.

MSC: 20M20; 20M18.
\end{abstract}

\section{Introduction}\label{sect:intro}

Fix a positive integer $n$, write $\bn=\{1,\ldots,n\}$, and denote by $\B_n$ the set of all binary relations on~$\bn$.  For $\al\in\B_n$ and for $x\in\bn$, write $x\al=\set{y\in\bn}{(x,y)\in\al}$ and $\al x=\set{y\in\bn}{(y,x)\in\al}$.  The set $\B_n$ forms a semigroup under the composition operation $\circ$ defined
by $\al\circ\be=\set{(x,y)\in\bn\times\bn}{x\al\cap\be y\not=\emptyset}$.  
In \cite{Vernitski2007}, the second author introduced and studied an alternative operation $\dia$ on $\B_n$, defined
by
\[
\al\dia\be=\set{(x,y)\in\bn\times\bn}{x\al=\be y\not=\emptyset}.
\]
It was shown in \cite{Vernitski2007} that the operation $\dia$ is not associative on $\B_n$, but that it is associative on the subset~$\D_n$ of~$\B_n$ consisting of all \emph{difunctional} relations on~$\bn$; see Section \ref{sect:prelim} for the definition of difunctionality.  
The semigroup $(\D_n,\dia)$ was shown to be an inverse semigroup in \cite{Vernitski2007}, and further properties of this semigroup were investigated in \cite{KM2011}, including Green's relations, ideals, maximal subsemigroups and congruences.  It was left as an open problem in \cite{KM2011} to determine the \emph{rank} of $\D_n$: that is, the minimal size of a (semigroup) generating set for~$\D_n$.\footnote{We note that Proposition 7 in an earlier version of \cite{KM2011}, available at {\tt arxiv.org/pdf/math/0602623v1.pdf}, leads to a lower bound for $\rank(\D_n)$ that is fairly close to the precise value.}  In this note, we solve this problem; see Theorem \ref{thm:main}.  In fact, we solve a more general problem, and calculate the rank of each ideal of $\D_n$; see Proposition \ref{prop:main}.  This being trivial for $n=1$, we assume $n\geq2$ for the remainder of the article.


\section{Preliminaries and statement of the main results}\label{sect:prelim}

Recall from \cite{Riguet1948} that a relation $\al$ on $\bn$ is \emph{difunctional} if $\al=\al\circ\al^{-1}\circ\al$, where $\al^{-1}=\set{(y,x)}{(x,y)\in\al}$ is the inverse relation of $\al$.
There are many equivalent formulations of the difunctionality property.  To describe the one that is most convenient for our purposes, we first introduce some notation.
For a set $X$, we write $\Part(X)$ for the set of all set partitions of $X$.  For $1\leq k\leq|X|$, we write $\Part(X,k)$ for the set of all set partitions of $X$ into $k$ blocks.  By convention, we also define $\Part(\emptyset)=\Part(\emptyset,0)=\{\emptyset\}$.  
%

A binary relation $\al\in\B_n$ is difunctional if and only if it is of the form $\al=(A_1\times B_1)\cup\cdots\cup(A_r\times B_r)$, for some subsets $A,B\sub\bn$ and some partitions $\{A_1,\ldots,A_r\}\in\Part(A,r)$ and $\{B_1,\ldots,B_r\}\in\Part(B,r)$.  We denote $\al$ as above by $\direl{A_1}{\cdots}{A_r}{B_1}{\cdots}{B_r}$.  We write
\begin{align*}
\rank(\al) &=r ,&
\dom(\al)&=A_1\cup\cdots\cup A_r ,&
\ker(\al)&=\{A_1,\ldots,A_r\} ,&
\deff(\al)&=|\bn\sm\dom(\al)|, \\
&& \codom(\al)&=B_1\cup\cdots\cup B_r ,&
\coker(\al)&=\{B_1,\ldots,B_r\} ,&
\codef(\al)&=|\bn\sm\codom(\al)|,
\end{align*}
and we call these parameters the \emph{rank}, \emph{domain}, \emph{codomain}, \emph{kernel}, \emph{cokernel}, \emph{defect} and \emph{codefect} of $\al$, respectively.  
Note that the empty relation $\emptyset$ is difunctional, corresponding to the $r=0$ case above.

Denote by $\I_n$ the subset of $\D_n$ consisting of all difunctional relations $\direl{A_1}{\cdots}{A_r}{B_1}{\cdots}{B_r}$ for which $|A_i|=|B_i|=1$ for each $1\leq i\leq r$.  It was shown in \cite{Vernitski2007} that $(\I_n,\dia)$ is a subsemigroup of $(\D_n,\dia)$; in fact, it was shown that the operations $\dia$ and $\circ$ coincide on $\I_n$, so that $\I_n$ is precisely the \emph{symmetric inverse monoid} on $\bn$.
In particular, the \emph{symmetric group} $\S_n=\set{\al\in\D_n}{\rank(\al)=n}$ is contained in $\D_n$.  
We note that the identity element of $\S_n$ is not an identity element of $\D_n$.  
In fact, $\D_n$ does not have an identity element,
but it does have a zero element, namely the empty relation, $\emptyset$.

Let $S$ be a semigroup, and write $S^1$ for the monoid obtained by adjoining an identity element to $S$ if necessary.  Recall that \emph{Green's preorders} $\leq_\R$, $\leq_\L$, $\leq_\J$ are defined, for $a,b\in S$ by
\[
a\leq_\R b \iff a\in bS^1 \COMMA
a\leq_\L b \iff a\in S^1b \COMMA
a\leq_\J b \iff a\in S^1bS^1,
\]
and that Green's relations $\R$, $\L$, $\J$ are defined by ${\R}={\leq_\R}\cap{\geq_\R}$, ${\L}={\leq_\L}\cap{\geq_\L}$, ${\J}={\leq_\J}\cap{\geq_\J}$.  Green's relation $\H$ is defined by ${\H}={\R}\cap{\L}$.  For more on Green's relations, and (inverse) semigroups more generally, the reader is referred to \cite{Howie1995,Lawson1998}.  The next result describes Green's relations and preorders on $\D_n$; its proof is routine, and is ommitted.  (Parts (iv)--(vi) may be found in \cite{KM2011}, in slightly different language, without proof.)

\ms
\begin{lemma}\label{lem:green}
Let $\al,\be\in\D_n$.  Then
\begin{itemize}\begin{multicols}2
\itemit{i} $\al\leq_\R\be$ if and only if $\ker(\al)\sub\ker(\be)$,
\itemit{ii} $\al\leq_\L\be$ if and only if $\coker(\al)\sub\coker(\be)$,
\itemit{iii} $\al\leq_\J\be$ if and only if $\rank(\al)\leq\rank(\be)$,
\itemit{iv} $\al\R\be$ if and only if $\ker(\al)=\ker(\be)$,
\itemit{v} $\al\L\be$ if and only if $\coker(\al)=\coker(\be)$,
\itemit{vi} $\al\J\be$ if and only if $\rank(\al)=\rank(\be)$.  \epfres
\end{multicols}\eit
\end{lemma}

It follows from parts (iii) and (vi) of Lemma \ref{lem:green} that the $\J$-classes of $\D_n$ are  the sets
\begin{align*}
J_r &= \set{\al\in\D_n}{\rank(\al)=r} &&\text{for $0\leq r\leq n$,}
\intertext{and that these form a chain under the usual ordering on $\J$-classes: $J_0<J_1<\cdots<J_n$.  That is, $J_r\sub \D_n\dia J_s\dia\D_n$ for any $0\leq r\leq s\leq n$.  Note also that $J_n=\S_n$ and $J_0=\{\emptyset\}$.  In any semigroup in which the $\J$-classes form a chain, the ideals form a chain under inclusion.  So the ideals of $\D_n$ are  the sets}
I_r = J_0\cup\cdots\cup J_r &= \set{\al\in\D_n}{\rank(\al)\leq r} &&\text{for $0\leq r\leq n$.}
\end{align*}
Our main results calculate the \emph{ranks} of these ideals, including that of $I_n=\D_n$ itself.  Recall that the rank of a semigroup $S$ is defined to be $\rank(S)=\min\bigset{|A|}{A\sub S,\ S=\la A\ra}$, the least cardinality of a generating set for $S$.  The rank of a semigroup should not be confused with the rank of a difunctional relation.  

To state our main results, we recall the definition of the Stirling and Bell numbers.
For non-negative integers $n$ and $k$, the \emph{Stirling number of the second kind} $S(n,k)$ denotes the number of partitions of a set of size $n$ into $k$ (nonempty) subsets.  The \emph{Bell number} $B(n)=S(n,1)+\cdots+S(n,n)$ denotes the total number of partitions of a set of size $n$ into any number of subsets.  Note that $S(0,0)=1$, and $S(n,k)=0$ if $k>n$.  The Stirling and Bell numbers are listed as Sequences A008277 and A000110, respectively, on  \cite{OEIS}.



\ms
\begin{prop}\label{prop:main}
Let $n\geq2$ and $0\leq r\leq n$.  Then the rank of the ideal $I_r = \set{\al\in\D_n}{\rank(\al)\leq r}$ of $\D_n$ is given by
\[
\rank(I_r) = \begin{cases}
\rho_{nr} &\text{if $r=0$ or $r\geq3$}\\
\rho_{nr}-1 &\text{if $1\leq r\leq2$,}
\end{cases}
\]
where $\displaystyle{\rho_{nr}=r+(r+1)S(n,r+1)+\sum_{k=1}^rS(n,k)}$.
\end{prop}

Proposition \ref{prop:main} yields a formula for the rank of $\D_n$ itself, upon putting $r=n$.  This formula may be simplified, noting that $S(n,n+1)=0$, and that $\sum_{k=1}^nS(n,k)=B(n)$:


\ms
\begin{thm}\label{thm:main}
If $n\geq3$, then $\rank(\D_n)=B(n)+n$.  \epfres
\end{thm}

For completeness, we note that $\rank(\D_2)=B(2)+1=3$.  We end this section with a simple combinatorial lemma.

\ms
\begin{lemma}\label{lem:green2}
Let $0\leq r\leq n$.  Then 
\bit
\itemit{i} $J_r$ contains $(r+1)S(n,r+1)+S(n,r)$ $\R$-classes (and the same number of $\L$-classes), and
\itemit{ii} the $\H$-class of any idempotent from $J_r$ is isomorphic to the symmetric group $\S_r$.
\eit
\end{lemma}

\pf By Lemma \ref{lem:green}(iv), an $\R$-class in $J_r$ is uniquely determined by the kernel of each of its elements.  This is a partition $\bA=\{A_1,\ldots,A_r\}$ of some subset $A$ of $\bn$ for which $|A|\geq r$.  The number of such partitions with $|A|=n$ is equal to $S(n,r)$.  The number of such partitions with $|A|<n$ is $(r+1)S(n,r+1)$; indeed, to specify such a partition, we first partition $\bn$ into $r+1$ blocks and choose one of these not to include as a block of~$\bA$.  This completes the proof of (i).

By Lemma \ref{lem:green}(iv) and (v), it is clear that the $\H$-class of the idempotent $\direl1\cdots r1\cdots r\in J_r$ consists of all permutations of the set $\{1,\ldots,r\}$,
so that part (ii) of the current lemma is true of this idempotent.  But all group $\H$-classes in $J_r$ are isomorphic; see \cite[Proposition 2.3.6]{Howie1995}. \epf

\ms
\begin{rem}
An alternative way of counting the $\R$-classes in $J_r$ involves (in the notation of the proof of Lemma \ref{lem:green2}) first choosing the subset $A$ and then the partition $\bA\in\Part(A,r)$.  This leads to the alternative expression of $\sum_{k=r}^n\binom nkS(k,r)$ for the number of such $\R$-classes.
\end{rem}

\section{Proof of the main result}\label{sect:proof}

Note that Proposition \ref{prop:main} is trivial for $r=0$, since $I_0=\{\emptyset\}$ and $\rho_{n0}=1$, so for the duration of this section, we fix $n\geq2$ and some $1\leq r\leq n$.  

Recall from \cite[Section 3.1]{Howie1995} that the \emph{principal factor} of a $\J$-class $J$ in a semigroup $S$ is the semigroup $J^*$ with underlying set $J\cup\{0\}$, where $0$ is a symbol not in $J$, and with product $*$ defined by
\[
a*b = \begin{cases}
ab &\text{if $a,b,ab\in J$}\\
0 &\text{otherwise.}
\end{cases}
\]
Recall from \cite{HRH1998} that the \emph{relative rank} of a semigroup $S$ with respect to a subset $A\sub S$, denoted $\rank(S:A)$, is the smallest cardinality of a subset $B\sub S$ such that $S=\la A\cup B\ra$.  The proof of the next result is routine, but is included for convenience.

\ms
\begin{lemma}\label{lem:SJ}
Let $S$ be a finite semigroup with a single maximal $\J$-class $J$ that is not a subsemigroup of~$S$.  Then $\rank(S)=\rank(J^*)+\rank(S:J)$.
\end{lemma}

\pf To avoid confusion during the proof, if $X\sub J$, we will write $\la X\ra$ for the subsemigroup of $S$ generated by $X$, and $\la X\ra^*$ for the subsemigroup of $J^*$ generated by $X$.

First, suppose $J^*=\la A\ra^*$ and $S=\la J\cup B\ra$, with $|A|=\rank(J^*)$ and $|B|=\rank(S:J)$.  Since $J$ is not a subsemigroup of $S$, we have $A\sub J$.
Then $\la A\cup B\ra = 
\la \la A\ra\cup B\ra=
\la J\cup B\ra=S$, so that $\rank(S)\leq|A\cup B|\leq|A|+|B|=\rank(J^*)+\rank(S:J)$.

Conversely, suppose $S=\la C\ra$, and put $A=C\cap J$ and $B=C\sm J$.  Let $x\in J$, and consider an expression $x=c_1\cdots c_k$, where $c_1,\ldots,c_k\in C$.  Since $S\sm J$ is a (nonempty) ideal of $S$, each factor $c_i$ must belong to $J$; that is $c_i\in A$.  It follows that $J\sub\la A\ra$, and so $J^*=\la A\ra^*$; note that $0\in \la A\ra^*$ because $J$ is not a subsemigroup of $S$.  In particular, $|A|\geq\rank(J^*)$.  But also $S=\la A\cup B\ra=\la\la A\ra\cup B\ra\supseteq\la J\cup B\ra\supseteq S$, so it follows that $S=\la J\cup B\ra$, giving $|B|\geq\rank(S:J)$.  Thus, $|C|=|A|+|B|\geq\rank(J^*)+\rank(S:J)$.  Since this is true for any generating set $C$ for $S$, it follows that $\rank(S)\geq\rank(J^*)+\rank(S:J)$. \epf

%

In the case that $S$ is the ideal $I_r$ of $\D_n$, it follows that $\rank(I_r)=\rank(J_r^*)+\rank(I_r:J_r)$.  We give the values of $\rank(J_r^*)$ and $\rank(I_r:J_r)$ in Lemmas \ref{lem:J} and \ref{lem:IJ}, respectively.

\ms
\begin{lemma}\label{lem:J}
If $1\leq r\leq n$, then $\rank(J_r^*)=\rank(\S_r)-1+(r+1)S(n,r+1)+S(n,r)$.
\end{lemma}

\pf Since $\D_n$ is an inverse semigroup, $J_r^*$ is a \emph{Brandt semigroup}.  More specifically, by Lemma \ref{lem:green2}(ii), $J_r^*$ is a Brandt semigroup over the symmetric group $\S_r$.  By \cite[Corollary 9]{Gray2014}, it follows that $\rank(J_r^*)=\rank(\S_r)-1+q$, where $q$ is the number of $\R$-classes in $J_r$.  The result now follows from Lemma \ref{lem:green2}(i). \epf

In light of Lemmas \ref{lem:SJ} and \ref{lem:J}, and the fact \cite{Moore1897} that
\[
\rank(\S_r) = \begin{cases}
1 &\text{if $r\leq 2$}\\
2 &\text{if $r\geq3$,}
\end{cases}
\]
the proof of Proposition \ref{prop:main} will be complete if we can prove the following.

\ms
\begin{lemma}\label{lem:IJ}
If $1\leq r\leq n$, then $\rank(I_r:J_r)=r-1+\displaystyle{\sum_{k=1}^{r-1}S(n,k)}$.
\end{lemma}


To prove Lemma \ref{lem:IJ}, we will first need to prove a number of intermediate results.  Consider a partition $\bA=\{A_1,\ldots,A_k\}\in\Part(\bn)$ with $\min(A_1)<\cdots<\min(A_k)$.  We define the difunctional relations
\[
\lam_\bA = \direl{A_1}\cdots{A_k}1\cdots k \AND \rho_\bA = \direl1\cdots k{A_1}\cdots{A_k}.
\]
Here and elsewhere, we use an obvious shorthand notation: for example, $\direl{A_1}\cdots{A_k}1\cdots k$ is an abbreviation for $\direl{A_1}\cdots{A_k}{\{1\}}\cdots{\{k\}}$.   For $1\leq k\leq n$, put
\[
\cL_k=\set{\lam_\bA}{\bA\in\Part(\bn),\ |\bA|\leq k} \AND \cR_k=\set{\rho_\bA}{\bA\in\Part(\bn),\ |\bA|\leq k}.
\]
Recall that the symmetric inverse monoid $\I_n$ is a subsemigroup of $\D_n$.  

\ms
\begin{lemma}\label{lem:LIR}
Let $\al\in I_{r-1}$.  Then $\al=\be\dia\ga\dia\de$ for some $\be\in\cL_r$, $\ga\in\I_n$, $\de\in\cR_r$ with ${\rank(\ga)=\rank(\al)}$.
\end{lemma}

\pf Write $\al=\direl{A_1}{\cdots}{A_k}{B_1}{\cdots}{B_k}$, noting that $k\leq r-1$.
Put $A_{k+1}=\bn\sm\dom(\al)$ and $B_{k+1}=\bn\sm\codom(\al)$, 
and let
\[
\bA = \begin{cases}
\{A_1,\ldots,A_k\} &\text{if $A_{k+1}=\emptyset$} \\
\{A_1,\ldots,A_k,A_{k+1}\} &\text{if $A_{k+1}\not=\emptyset$} 
\end{cases}
\AND
\bB = \begin{cases}
\{B_1,\ldots,B_k\} &\text{if $B_{k+1}=\emptyset$} \\
\{B_1,\ldots,B_k,B_{k+1}\} &\text{if $B_{k+1}\not=\emptyset$.} 
\end{cases}
\]
Then it is easy to see that $\al=\lam_\bA\dia\ga\dia\rho_\bB$, where $\ga=\rho_\bA\dia\al\dia\lam_\bB\in\I_n$ with $\rank(\ga)=k$. \epf


\ms
\begin{lemma}\label{lem:LSiR}
We have $I_r=\la J_r\cup\cL_r\cup\cR_r\ra$.
\end{lemma}

\pf We must consider two separate cases.  Suppose first that $r<n$.  Note that $J_r$ contains the set $\Om=\set{\al\in\I_n}{\rank(\al)=r}$.  It is well known that $\la\Om\ra=\set{\al\in\I_n}{\rank(\al)\leq r}$; see for example \cite[Lemma~4.7]{ZF2015}.  The result now follows from Lemma~\ref{lem:LIR}. 

Suppose now that $r=n$, so $J_r=\S_n$.  By Lemma \ref{lem:LIR}, it suffices to show that $\I_n\sub\la \S_n\cup\cL_n\cup\cR_n\ra$.  For this, let $\bA\in\Part(\bn,n-1)$ be arbitrary, and put $\al=\rho_\bA\dia\lam_\bA\in\la \cL_n\cup\cR_n\ra$, noting that $\al=\direl1\cdots{n-1}1\cdots{n-1}\in\I_n$ and $\rank(\al)=n-1$.  It then follows from the proof of \cite[Theorem 3.1]{GH1987} (see also \cite{Popova1961}) that $\I_n=\la\S_n\cup\{\al\}\ra$. \epf


Let $\bA,\bB\in\Part(\bn,r)$, and write $\bA=\{A_1,\ldots,A_r\}$ and $\bB=\{B_1,\ldots,B_r\}$ with $\min(A_1)<\cdots<\min(A_r)$ and $\min(B_1)<\cdots<\min(B_r)$.  We define $\phi_{\bA,\bB}=\direl{A_1}{\cdots}{A_r}{B_1}{\cdots}{B_r}$.  In order to simplify notation in what follows, and since $n$ is fixed, for each $1\leq k\leq n$, we will write $p_k=S(n,k)$.  For each $1\leq k\leq n$, let us denote the elements of $\Part(\bn,k)$ by $\bA_{k,1},\ldots,\bA_{k,p_k}$.

\ms
\begin{lemma}\label{lem:upper}
For each $1\leq k\leq n-1$, let $\Si_k = \{\phi_{\bA_{k,1},\bA_{k,2}},\ldots,\phi_{\bA_{k,p_k-1},\bA_{k,p_k}}\} \cup \{\lam_{\bA_{k,p_k}},\rho_{\bA_{k,1}}\}$.
Then for any $1\leq r\leq n$,
$I_r=\la J_r\cup\Si_1\cup\cdots\cup\Si_{r-1}\ra$.
\end{lemma}

\pf Put $\Om = J_r\cup\Si_1\cup\cdots\cup\Si_{r-1}$.  By Lemma \ref{lem:LSiR}, to show that $I_r=\la\Om\ra$, it suffices to show that~$\la\Om\ra$ contains both $\cL_r$ and $\cR_r$.  Let $\bA\in\Part(\bn)$ with $|\bA|\leq r$.  We must show that $\lam_\bA,\rho_\bA\in\la\Om\ra$.  If $|\bA|=r$, then $\lam_\bA,\rho_\bA\in J_r\sub\la\Om\ra$, so suppose $\bA\in\Part(\bn,k)$, where $1\leq k\leq r-1$.  Then $\bA=\bA_{k,l}$ for some $1\leq l\leq p_k$.  But then
\begin{align*}
\lam_\bA = \lam_{\bA_{k,l}} &= (\phi_{\bA_{k,l},\bA_{k,l+1}} \dia\cdots\dia \phi_{\bA_{k,p_k-1},\bA_{k,p_k}}) \dia \lam_{\bA_{k,p_k}} , \\
\rho_\bA = \rho_{\bA_{k,l}} &= \rho_{\bA_{k,1}} \dia (\phi_{\bA_{k,1},\bA_{k,2}} \dia\cdots\dia \phi_{\bA_{k,l-1},\bA_{k,l}}) ,
\end{align*}
where the first bracketed expression is omitted if $l=p_k$, and the second if $l=1$. \epf

\ms
\begin{rem}
We could not help noticing that the generating set used in Lemma \ref{lem:upper} looks very similar to the construction of so-called \emph{rainbow tables} in computer security \cite{Hellman1980}. This is perhaps not surprising, since both constructions have the purpose, broadly speaking, of reducing the total amount of memory used for storing given information.
\end{rem}

Since $|\Si_k|=S(n,k)+1$ for each $k$, it follows from Lemma \ref{lem:upper} that $\rank(I_r:J_r)\leq r-1+\sum_{k=1}^{r-1}S(n,k)$.  To complete the proof of Lemma \ref{lem:IJ}, we must therefore show that this upper bound for $\rank(I_r:J_r)$ is also a lower bound.  To do this, we will show in Lemmas \ref{lem:Si1} and \ref{lem:Si2} that if $\Si\sub I_r$ is such that $I_r=\la J_r\cup\Si\ra$, then $\Si$ must include certain specified kinds of relations.  
First, we prove two intermediate lemmas.  There are obvious dual versions of Lemmas \ref{lem:abc} and \ref{lem:abc2}, but we will not state them.

\ms
\begin{lemma}\label{lem:abc}
If $\al,\be,\ga\in\D_n$ are such that $\al=\be\dia\ga$ and $\dom(\al)=\bn$, then $\ker(\al)=\ker(\be)$.
\end{lemma}

\pf Since $\al=\be\dia\ga$, we have $\al\leq_\R\be$, so Lemma \ref{lem:green}(i) gives $\ker(\al)\sub\ker(\be)$.  Since $\dom(\al)=\bn$, it is clear that $\ker(\al)$ is maximal, inclusion-wise, so we must in fact have $\ker(\al)=\ker(\be)$. \epf

%

\ms
\begin{lemma}\label{lem:abc2}
If $\al,\be,\ga\in\D_n$ are such that $\al=\be\dia\ga$, $\ker(\al)=\ker(\be)$ and $\codom(\be)=\bn$, then 
$\be^{-1}\dia\al=\ga$.
\end{lemma}

\pf Since $\ker(\be)=\ker(\al)=\ker(\be\dia\ga)$, it follows that $\coker(\be)\sub\ker(\ga)$.  Since $\codom(\be)=\bn$, $\coker(\be)$ is maximal, inclusion-wise, so we must in fact have $\coker(\be)=\ker(\ga)$.  But then $\be^{-1}\dia\be=\ga\dia\ga^{-1}$, which gives $\ga=\ga\dia\ga^{-1}\dia\ga=\be^{-1}\dia\be\dia\ga=\be^{-1}\dia\al$. \epf


\ms
\begin{lemma}\label{lem:Si1}
If $I_r=\la J_r\cup\Si\ra$, and if $1\leq k\leq r-1$, then there exist $\si,\tau\in\Si$ with $\dom(\si)=\codom(\tau)=\bn$, $\rank(\si)=\rank(\tau)=k$ and $\codef(\si),\deff(\tau)>0$.
\end{lemma}

\pf It suffices to prove the existence of $\si$, as the existence of $\tau$ will follow by a symmetrical argument (for which we need the duals of Lemmas \ref{lem:abc} and \ref{lem:abc2}).  Let $1\leq k\leq r-1$, and write 
\[
\Om=\set{\al\in\D_n}{\dom(\al)=\bn,\ \rank(\al)=k,\ \codef(\al)>0}.
\]
For $\al\in\Om$, write $\ell(\al)$ for the minimum value of $m$ such that $\al=\be_1\dia\cdots\dia\be_m$ for some $\be_1,\ldots,\be_m\in J_r\cup\Si$.  Let $L=\min\set{\ell(\al)}{\al\in\Om}$.  
To establish the existence of $\si$, it suffices to prove that $L=1$.  
%
To do this, suppose to the contrary that $L\geq2$, and choose some $\al= \direl{A_1}\cdots{A_k}{B_1}\cdots{B_k}\in\Om$ with $\ell(\al)=L$.  So we may write $\al=\be_1\dia\be_2\dia\cdots\dia\be_L$ for some $\be_1,\be_2,\ldots,\be_L\in J_r\cup\Si$.   For simplicity, put $\be=\be_1$ and $\ga=\be_2\dia\cdots\dia\be_L$, so $\al=\be\dia\ga$.  Lemma \ref{lem:abc} gives $\ker(\be)=\ker(\al)$, so we may write $\be=\direl{A_1}\cdots{A_k}{C_1}\cdots{C_k}$.  If $\codef(\be)>0$, then we put $\si=\be$, and the proof of the lemma is complete.  So suppose $\codef(\be)=0$.  This means that $\codom(\be)=\bn$, and Lemma \ref{lem:abc2} then gives 
\begin{equation}
\label{eq:*}
\be^{-1}\dia\al=\ga=\be_2\dia\cdots\dia\be_L.
\end{equation}
But $\be^{-1}\dia\al=\direl{C_1}\cdots{C_k}{A_1}\cdots{A_k}\dia\direl{A_1}\cdots{A_k}{B_1}\cdots{B_k}=\direl{C_1}\cdots{C_k}{B_1}\cdots{B_k}$.  
Consequently,  $\dom(\be^{-1}\dia\al)=
\codom(\be)=\bn$ and $\codef(\be^{-1}\dia\al)=
\codef(\al)>0$.  Thus, $\be^{-1}\dia\al\in\Om$.  But $\ell(\be^{-1}\dia\al)\leq L-1$, by \eqref{eq:*}, contradicting the minimality of $L$.  This completes the proof. \epf

\ms
\begin{lemma}\label{lem:Si2}
If $I_r=\la J_r\cup\Si\ra$, and if $\bA\in\Part(\bn)$ with $|\bA|\leq r-1$, then there exist $\si,\tau\in\Si$ with $\ker(\si)=\bA$ and $\coker(\tau)=\bA$.
\end{lemma}

\pf Again, it suffices to demonstrate the existence of $\si$.  Choose some $\al\in I_r$ with $\ker(\al)=\bA$, noting that $\dom(\al)=\bn$.  Suppose $\al=\be_1\dia\cdots\dia\be_k$ where $\be_1,\ldots,\be_k\in J_r\cup\Si$.  If $k=1$, then $\al=\be_1\in\Si$, and we are done, with $\si=\al$.  If $k\geq2$, then $\al=\be_1\dia(\be_2\dia\cdots\dia\be_k)$, and Lemma \ref{lem:abc} gives $\ker(\be_1)=\ker(\al)=\bA$, and we are done with $\si=\be_1$. \epf

\pfof{Lemma \ref{lem:IJ}} 
As noted after the proof of Lemma \ref{lem:upper}, it suffices to show that $\rank(I_r:J_r)\geq r-1+\sum_{k=1}^{r-1}S(n,k)$.  Suppose $I_r=\la J_r\cup\Si\ra$.  
For each $1\leq k\leq r-1$, let $\Si_k=\set{\al\in\Si}{\rank(\al)=k}$, and fix some such $k$.  It is enough to show that $|\Si_k|\geq1+S(n,k)$.  By Lemma \ref{lem:Si1}, there exists some $\tau\in\Si_k$ with $\deff(\tau)>0$.  By Lemma \ref{lem:Si2}, for any $\bA\in\Part(\bn,k)$, there exists some $\si_\bA\in\Si_k$ with $\ker(\si_\bA)=\bA$.  Clearly these elements of $\Si$ are all distinct, so $|\Si_k|\geq1+|\Part(\bn,k)|=1+S(n,k)$, as required. \epf

As noted before the statement of Lemma \ref{lem:IJ}, this completes the proof of Proposition \ref{prop:main}.

\ms
\begin{rem}
%
Finally, we note that $\D_n$ bears many similarities with several families of semigroups, such as the symmetric inverse monoids $\I_n$, the full and partial transformation monoids $\T_n$ and $\P\T_n$, and certain diagram monoids such as the partition monoids $\P_n$.  All these monoids have a chain of $\J$-classes, and have the symmetric group~$\S_n$ as their (unique) maximal $\J$-class.  However, the ranks of the monoids $\I_n$, $\T_n$, $\P\T_n$ and $\P_n$ are constant and very small (all being equal to either $3$ or $4$, for $n\geq3$), and each monoid may be generated by elements in its top two $\J$-classes; see \cite{Aizenstat1958,GH1987,JEgrpm,Vorobev1953}.  The proper ideals of these monoids are all generated by elements in a single $\J$-class; formulae for the ranks of the ideals of these monoids may be found in \cite{EG2017,ZF2015,HM1990,Garba1990}.  By contrast, as we have seen, $\rank(\D_n)=B(n)+n$ grows rapidly with $n$, and any generating set for $\D_n$ or one of its proper ideals must contain elements from all $\J$-classes except the very bottom one.
Calculated values of $\rank(I_r)$ and $\rank(\D_n)$ are given in Tables \ref{tab:rankIr} and \ref{tab:rankDn}, respectively.
\end{rem}

\begin{table}[h]%
\begin{center}
\begin{tabular}{|c|rrrrrrrrrrr|}
\hline
$n\setminus r$&0&1&2&3&4&5&6&7&8&9&10 \\ 
\hline
\phantom{1}0& 1 &&&&&&&&&&\\
\phantom{1}1& 1& 2 &&&&&&&&&\\
\phantom{1}2& 1& 3& 3 &&&&&&&&\\
\phantom{1}3& 1& 7& 8& 8 &&&&&&&\\
\phantom{1}4& 1& 15& 27& 21& 19 &&&&&&\\
\phantom{1}5& 1& 31& 92& 84& 60& 57 &&&&&\\
\phantom{1}6& 1& 63& 303& 385& 266& 213& 209 &&&&\\
\phantom{1}7& 1& 127& 968& 1768& 1419& 986& 889& 884 &&&\\
\phantom{1}8& 1& 255& 3027& 7901& 8049& 5446& 4313& 4154& 4148 &&\\
\phantom{1}9& 1& 511& 9332& 34364& 45810& 33883& 23888& 21405& 21163& 21156 &\\
10& 1& 1023& 28503& 146265& 256576& 223439& 150465& 121186& 116342& 115993& 115985 \\
\hline
\end{tabular}
\end{center}
\vspace{-5mm}
\caption{Values of $\rank(I_r)$; see Proposition \ref{prop:main}.}
\label{tab:rankIr}
\end{table}

\begin{table}[h]%
\begin{center}
\begin{tabular}{|c|rrrrrrrrrrrrrr|}
\hline
$n$&0&1&2&3&4&5&6&7&8&9&10 &11&12&13 \\ 
\hline
$\rank(\D_n)$ &  1& 2& 3& 8& 19& 57& 209& 884& 4148& 21156& 115985& 678581& 4213609& 27644450 \\
\hline
\end{tabular}
\end{center}
\vspace{-5mm}
\caption{Values of $\rank(\D_n)$; see Theorem \ref{thm:main}.}
\label{tab:rankDn}
\end{table}

\footnotesize
\def\bibspacing{-1.1pt}
\bibliography{biblio}
\bibliographystyle{plain}
\end{document}